\documentstyle[12pt]{article}
\newtheorem{theorem}{THEOREM}[section]

\def\wbar{\overline{w}}

\def\dist{\hbox{dist}}

\def\hbar{\overline{h}}

\def\Omegabar{\overline{\Omega}}

\font\teneufm=eufm10
\font\seveneufm=eufm7
\font\fiveeufm=eufm5
\newfam\eufmfam
\textfont\eufmfam=\teneufm
\scriptfont\eufmfam=\seveneufm
\scriptscriptfont\eufmfam=\fiveeufm

\newfam\msbfam
\font\tenmsb=msbm10 scaled \magstep1   \textfont\msbfam=\tenmsb
\font\sevenmsb=msbm7 scaled \magstep1  \scriptfont\msbfam=\sevenmsb
\font\fivemsb=msbm5 scaled \magstep1   \scriptscriptfont\msbfam=\fivemsb
\def\Bbb{\fam\msbfam \tenmsb}

\def\RR{{\Bbb R}}
\def\CC{{\Bbb C}}

\def\PP{{\Bbb P}}

\def\hexdigit#1{\ifnum#1<10 \number#1\else
 \ifnum#1=10 A\else\ifnum#1=11 B\else\ifnum#1=12 C\else
 \ifnum#1=13 D\else\ifnum#1=14 E\else\ifnum#1=15 F\fi\fi\fi\fi\fi\fi\fi}

\chardef\\="5C                    %% Typesets \ in \tt mode
\chardef\{="7B  \chardef\}="7D    %% also left and right braces

 \def\HollowBoxx #1#2#3{{\dimen0=#1 \advance\dimen0 by -#2       
       \dimen1=#1 \advance\dimen1 by #3                       
        \vrule height 0pt depth #3 width #2                   
       \hskip -#3
       \vrule height #1 depth #3 width #3}}                   
 \def\LeftContraction{\mathord{\kern1.45pt \HollowBoxx{6pt}{3.5pt}{.4pt}}\,
}

 \def\HollowBox #1#2#3{{\dimen0=#1 \advance\dimen0 by -#3       
       \dimen1=#1 \advance\dimen1 by #3                       
        \vrule height #1 depth #3 width #3                    
        \vrule height 0pt depth #3 width #2                   
        \hskip -#3}}                                             
 \def\RightContraction{\mathord{\, \HollowBox{6pt}{3.1pt}{.4pt}}
 \kern1.6pt}

\def\jbar{\overline{j}}

\def\hbar{\overline{h}}

\def\Omegabar{\overline{\Omega}}

\def\epf{\hskip.2in\vrule width.4pt height6.65pt
depth.15pt\vrule
width2.5pt height6.65pt depth-6.25pt\hskip-2.5pt\vrule
width2.5pt
height.25pt depth.15pt\vrule width.4pt
height6.65pt depth.15pt\ }

\def\epf{\hskip.2in\vrule width.4pt height6.65pt
depth.15pt\vrule
width2.5pt height6.65pt depth-6.25pt\hskip-2.5pt\vrule
width2.5pt
height.25pt depth.15pt\vrule width.4pt
height6.65pt depth.15pt\ }

\def\proof{\noindent {\bf Proof. }}

\def\nbar{\overline{n}}
\def \Omegabar{\overline \Omega}

\def \d{\partial}

\def\jbar{\overline{j}}

\def\wbar{\overline{w}}

\def\zbar{\overline{z}}
\def\dbar{\overline{\partial}}
\def \hbar{\overline{h}}

\def\jbar{\overline{j}}

\def\\GBB{\cal B}

\begin{document}

\begin{center}
\large \bf \boldmath Application of the Complex Monge-Amp\`{e}re Equation
to the Study of  Proper Holomorphic Mappings
of Strictly Pseudoconvex Domains 
\medskip \\
\normalsize Steven G. Krantz\footnote{Author partially
supported by NSF Grant DMS--9531967.
Work at MSRI supported by NSF Grant DMS--9022140.} 
%% Song-Ying:  Do you prefer to omit this reference to MSRI?
%% It's not really important to me, but if we omit it then
%% we can't post our paper on the MSRI preprint server, which
%% is actually a server that a lot of people look at. Just
%% let me know what you think. Just fine.
\ \ and \ \  Song-Ying  Li
\footnote{Author partially
supported by NSF Grant DMS--9500758.}
\smallskip \\
\small December 2, 1996
\end{center}

\begin{quote}
\small \bf Abstract:  \rm We construct a special plurisubharmonic
defining function for a smoothly bounded strictly pseudoconvex
domain so that the determinant of the complex Hessian
vanishes to high order on the boundary.  This construction,
coupled with regularity of solutions of complex Monge-Amp\' ere
equation and the reflection principle,
enables us to give a new proof of the Fefferman mapping theorem.
\end{quote}
\vspace*{.12in}

\section{Introduction} 

In classical analysis, an important theorem of 
Painlev\'{e} [Pai] and Kellogg [Kel] states
that any conformal mapping between two smoothly bounded domains in
the complex plane $\CC$ can be extended to be a diffeomorphism on the
closures of the domains. This theorem was generalized by C.\
Fefferman [Fef] in 1974 to strictly pseudoconvex domains in $\CC^n$.
Fefferman's original proof of this theorem is very technical,
relying as it does on deep work on
the boundary asymptotics of the Bergman
kernel and on the
regularity of $\dbar$-Neumann operator that is due to J.\ J.\ Kohn [Ko]. 
Bell/Ligocka [BeL], and later Bell [Be], gave a simpler
proof which deals with more general domains, 
including pseudoconvex domains of finite type,
by using regularity of the Bergman projection and the 
$\dbar$-Neumann operator
as studied by Kohn [Ko], Catlin [Ca], Boas/Straube [BS], and
others.

From the results in  [NWY], [Web], [Lem],  [PH] and [For], we know that
 the  proof of the Fefferman mapping
theorem can be reduced to proving that $\det(\varphi'(z))$ 
extends continuously
across the boundary if $\varphi:\Omega_1\to \Omega_2$ 
is a biholomorphic mapping and
the $\Omega_j, j = 1,2, $ are 
smoothly bounded strictly pseudoconvex domains in
 $\CC^n$;
such an argument uses the reflection principle as developed
by  Nirenberg/Webster/Yang [NWY], Webster [Web], 
Pinchuk/Hasanov [PH],
Coupet [Cou], and more recently F.\ Forstneric [For]. 

We know from [Ker] that Painlev\'{e} and Kellogg's
theorem can be proved by using the regularity of the Dirichlet
problem for the Laplacian in a smoothly bounded planar domain, where
the property of the Laplacian being conformally invariant plays an
important role in the proof. The natural generalization of the
Laplacian in one complex variable to several complex variables, with
these considerations in mind, is the complex Monge-Amp\`{e}re
equation. In [Ker], Kerzman observed that the proof of the Fefferman
mapping theorem would follow from
 the $C^{\infty}$ global regularity of the Dirichlet problem of a
degenerate  complex Monge-Amp\`ere equation. However, counterexamples
in Bedford/Forn\ae ss [BF] as well as in  Gamelin/Sibony [GS]
show that, in general, the degenerate Dirichlet problem for the complex
Monge-Amp\`ere equation does not have $C^2$ boundary regularity. Thus
Kerzman's idea does not work in the sense of its original
formulation.

 The main purpose of the present paper is to give a new method 
for using the complex Monge-Amp\`{e}re equation to study the
boundary regularity of biholomorphic mappings of strongly
pseudoconvex domains; thus, in effect, we validate the program
initiated by Kerzman/Kohn/Nirenberg [Ker].  We achieve
this goal by proving our Theorem 2.2 and then  combining our result
with   a result  of
 Caffarelli/Kohn/Nirenberg/Spruck [CKNS] on the
Dirichlet problem for the complex Monge-Amp\`{e}re equation.
Specifically, we shall prove the following theorem.

\begin{theorem} \sl Let $\Omega_j, j = 1,2$ be two bounded strictly
pseudoconvex domains in $\CC^n$ with $C^{\infty}$ boundary.
 Let
$\varphi :\Omega_1 \to \Omega_2$ be a proper holomorphic mapping.
Then  for any $\epsilon > 0$ we have
\smallskip 

(i) $\varphi$ can be extended as a $\hbox{Lip}_{1}(\Omegabar_1)$ 
mapping;
 
(ii) $\det(\varphi')\in \hbox{Lip}_{1/2}(\Omegabar_1)$;

(iii) There is a $C^2$ defining function $\rho$ for $\Omega_2$ so that,
for any $k = 1,
\dots, n$, we have $\sum_{j=1}^n \left [ {\d \rho \over \d w_j} \circ
\varphi \right ] \cdot {\d \varphi_j \over \d z_k} \in
\hbox{Lip}_{1}(\Omegabar_1)$.  \end{theorem}

As a corollary of (ii) and of a theorem in [PH] and [Web],
we obtain a new proof of the Fefferman mapping theorem. 

\medskip

This paper is organized as follows. In Section 2, our main theorem
(Theorem 2.2) is stated and proved. Theorem 1.1
 is proved in the second part of Section 2
 by combining Theorem 2.2 and Theorem 2.1
(a result in [CKNS]). We shall then include
an explanation (in Section 3) of how to combine theorems 
in [PH] and [Web] together with
Theorem 1.1 to obtain a new proof of the Fefferman mapping theorem.
\medskip

The authors would like to thank Xiaojun Huang
and Peter Li for helpful conversations.
              
\section{An Application of the Complex Monge-Amp\`ere Equation}

Let us recall a theorem of Caffarelli, Kohn, Nirenberg and Spruck [CKNS].

\begin{theorem}\sl Let $\Omega$ be a bounded, strictly pseudoconvex
domain in $\CC^n$ with $C^3$ boundary $\d \Omega$. Let $f$ be a
 non-negative
function on $\Omegabar$ such that $f(z)^{1/n}\in C^{1,1}(\Omegabar)$.
Let $H(u)$ denote the complex Hessian matrix of the function $u$.
Then there is a unique plurisubharmonic function 
$v\in C^{1,1}(\Omegabar)$ satisfying
\begin{eqnarray*}
\det H(v)(z)  & = & f(z) \quad \hbox{for}\ \ z \in \Omega; \\
           v  & = & 0 \ \quad \ \ \mbox{for} \ \  z \in \d \Omega.  
                             \qquad \qquad \qquad \qquad (2.1) 
\end{eqnarray*}
\end{theorem}

\noindent {\bf Note:} The theorem stated in [CKNS] 
is for $\d \Omega$ being $C^{3,1}$ but, we may
approximate $\d \Omega$ by a sequence of bounded strictly pseudoconvex
domains with $C^4$ boundaries, and do a careful count
of which derivatives are actually used in their proof,
to obtain the result stated here.
\smallskip 

As we mentioned in the introduction, the following 
theorem plays an essential role in the proof of Theorem 1.1.

\begin{theorem}\sl Let $\Omega$ be a bounded strictly pseudoconvex
domain in $\CC^n$ with $C^{\infty}$ boundary $\d \Omega$. For any 
$0<\epsilon << 1$
and any positive integer $q$ there is a 
plurisubharmonic defining function $\rho_q\in
C^{\infty}(\Omegabar)$ for $\Omega$ so that
$$
\det H(\rho_q)(z) \le C\dist(z,\d \Omega)^q,\quad z\in \Omega. \leqno(2.3)
$$
\end{theorem}

\proof 
Let $\delta(z)$
denote the (signed) distance function from $z$ to $\d \Omega$. 
Since $\Omega$ is a bounded domain in
$\CC^n$ with $C^{\infty}$ boundary $\d \Omega$, then
 $\delta(z)\in C^{\infty}(\Omegabar)$ (after modification
of the distance function on a compact set in the interior).  
By a rotation we see that,
for any fixed point $z^0\in \Omega$ near the boundary, we may assume that 
the $z_n$ direction
is the direction $({\d \delta(z)\over \d \zbar_1},\cdots,{\d \delta\over
\d \zbar_n})$ at the point $z^0$.  Let
$$
H(\delta)_{n-1}(z)=\left[{\d^2 \delta \over \d z_\ell \d \zbar_q}(z)
\right]_{1\le \ell,q\le n-1}
$$
Since $\Omega$ is strictly pseudoconvex, there is an $\epsilon > 0$
so that 
$$
-H(\delta)_{n-1}(z^0)\ge \epsilon I_{n-1}
$$
for all $z\in \Omegabar$ with $\delta(z)\le \epsilon$. 
We may assume that $H(\delta)_{n-1}$ is diagonal at $z^0$.
Note that
\begin{eqnarray*}
H(- \delta)(z)
&=& \left[\matrix{- H(\delta)_{n-1}(z) & 0\cr 0 & -{\d^2 \delta \over \d z_
n
\d \zbar_n }
\cr}\right]
\end{eqnarray*}
If it happens that
$\sum_{i j}{\d \delta\over \d z_i }{\d \delta \over \d \zbar_j}
 {\d^2 \delta \over \d z_i \d \zbar_j} = 0$,
then $\rho_q=-\delta(z)+ \delta(z)^{2+q}$ is the desired 
defining function.  [In fact it is these terms that
distinguish the study of the real Hessian from the
more subtle study of the complex Hessian.  In particular,
we know that ${\d^2 \delta\over \d x_n^2}$ equals $0$, but
the term ${\d^2 \delta \d z_n \d \zbar_n}$ may not be zero.
Therefore estimate (2.3) is easy to check
for the determinant of
the real Hessian of $\delta$; matters are much trickier
for the complex Hessian.]

Now we let
$$
\rho(z)=-\delta(z)
$$
and
$$
r^{[m]}(z)=\rho(z)+\sum_{k=2}^m a_k(z) \delta(z)^k \quad z\in \Omega.
$$
We will prove inductively that 
$$
\hbox{det} (H(r^{[m]})(z)) = b_{m-1}(z) \cdot \rho^{m-1} ,
$$
where $b_{m-1}$ is some smooth function.  In particular,
$\hbox{det} (H(r^{[m]}))$ vanishes to order $m-1$ at the boundary.

Now
\begin{eqnarray*}
\d_{i\jbar} r^{[m]}(z)
&=& \d_{i\jbar} \rho(z) + \sum_{k=2}^m
[\d_{i\jbar} a_k \, \rho(z)^2+k (k-1) a_k \d_i \rho
\d_{\jbar} \rho)] \rho^{k-2}\\
&& +\sum_{k=2}^m k[\d_i a_k \d_{\jbar} \rho
+\d_i \rho \d_{\jbar} a_k 
+ a_k\rho_{i\jbar}] \rho^{k-1}\\
\end{eqnarray*}
Let $H(r^{[m]})_{n-1}(z)=[{\d^2 r^{[m]} \over \d z_i \d
 \zbar_j}]_{(n-1)\times (n-1)}$,
and let $B(z)^*=[r_{n \bar{1}},\cdots, r_{n\overline{n-1}}]$ 
be a row-vector. Then
\begin{eqnarray*}
H(r^{[m]})(z)
&=& \left[\matrix{H(r^{[m]})_{n-1}(z) & B(z)\cr B(z)^* & r^{[m]}_{n\nbar}(z
) 
\cr}\right]
\end{eqnarray*}
Then
$$
\det( H(r^{[m]})(z))=\det (H(r^{[m]})_{n-1}(z))
[r^{[m]}_{n\nbar} - \{ B^* H(r^{[m]})_{n-1}(z)\}^{-1} B(z)]
$$
We know that
$$
H(r^{[m]})_{n-1}(z)\ge \epsilon I_{n-1}
$$
for all $z\in \Omega_{\epsilon}$, where $\epsilon$ is
a positive number depending only on $\Omega$ and 
$\|a_j\|_{C^2(\Omegabar)}$. First, we let
$$
(d\rho)(z)= \sum_{i,j}\rho_i \rho_{\jbar} \rho_{i\jbar} \, ;
$$
then we define 
$$
a_2(z)=-{(d \rho)(z)  \over 2
((d \rho)(z) +|\d \rho|^4)}.
$$
Thus we have
\begin{eqnarray*}
r_{n\nbar}(z)
&=&\rho_{n\nbar} +\sum_{k=2}^m
[ \d_{n \nbar} a_k \rho^k + k(k-1) a_k |\d \rho|^2 \rho^{k-2}]\\
&&+\sum_{k=2}^m k[\d_n a_k \d_{\nbar} \rho +
\d_n \rho \d_{\nbar} a_k +a_k \rho_{n\nbar}] \rho^{k-1}\\
&=&\sum_{k=2}^m
[ \d_{n \nbar} a_k \rho^k + k(k+1) a_{k+1} |\d \rho|^2 \rho^{k-1}]\\
&&+\sum_{k=2}^m k[\d_n a_k \d_{\nbar} \rho +
\d_n \rho \d_{\nbar} a_k +a_k \rho_{n\nbar}] \rho^{k-1}\\
&=& O(\rho)
\end{eqnarray*}
Therefore
$$
\det(H(r^{[2]})(z))=b_1(z) \rho(z)
$$
Assume that we have constructed $r^{[m]}=\rho(z)+\sum_{k=2}^m a_k
 \rho(z)^k$
such that
$$
\det (H(r^{[m]})(z))=b_{m-1}(z)\rho^{m-1},\quad z\in \Omega.
$$
We consider
$$
r^{[m+1]}=r^{[m]} +a_{m+1} \rho(z)^{m+1}.
$$
Since
\begin{eqnarray*}
\lefteqn{H(r^{[m+1]})(z)}\\
&=&H(r^{[m]})+\rho(z)^{m+1} H(a_{m+1})(z)
+(m+1) a_{m+1} \rho(z)^m H(\rho)\\
&& +(m+1) m a_{m+1}\rho^{m-1}
\d \rho(z) \otimes \overline{\d \rho} \\
&& +(m +1)\rho(z)^m [\d \rho \otimes \overline{
\d a_{m+1}}+\d a_{m+1} \otimes \overline{\d \rho}] ,
\end{eqnarray*}
it is easy to see that
$$
\det(H(r^{[m+1]})(z))
=b_m(z) \rho(z)^m)+\det
 \bigl [ H(r^{[m]})+ (m+1) m a_{m+1}\rho^{m-1}) \bigr ] .
$$
By a rotation, we may let $z_n$ be the complex normal direction of
 $\d \Omega_{\delta(z)}$ at $z$. Thus
\begin{eqnarray*}
\det(H(r^{[m]})(z)+ (m+1) m \rho(z)^{m-1}\d \rho(z)\otimes 
       \overline{\d \rho}) && \\
&=&\det\, M ,
\end{eqnarray*}
where
$$
M = \left [\matrix{ H(r^{[m]})_{n-1}(z) & B(m)(z)\cr
B(m)(z)^* & r^{[m]}_{n \nbar}(z)+ (m+1) m \rho^{m-1} \rho_n
 \rho_{\nbar}\cr} \right ] .
$$
But this 
\begin{eqnarray*}
&=& \det(H(r^{[m]})_{n-1}(z) [r^{[m]}_{n \nbar}+ B(m)^*
 H(r^{[m]})_{n-1}(z)^{-1}
B(m)(z) \\
&& \qquad +(m+1) m \rho^{m-1} \rho_n \rho_{\nbar}]  \\
&=& \det (H(r^{[m]})(z))+\det(H(r^{[m]})_{n-1}(z)) (m+1) m \rho_n
 \rho_{\nbar} \rho(z)^{m-1}\\
&=& b_{m-1}(z) \rho(z)^{m-1}+
\det(H(r^{[m]})_{n-1}(z)) (m+1) m \rho_n \rho_{\nbar} \rho(z)^{m-1} .
\end{eqnarray*}
We therefore choose $a_{m+1}$ such that
$$
b_{m-1}(z)+ (m+1) m \det(H(r^{[m]})_{n-1}(z)) \rho_n \rho_{\nbar}(z)=0
$$
From this it will follow that
$$
\det(H(r^{[m+1]})(z))=b_m (z) \rho(z)^m.
$$
By our construction, we know that $H(r^{[m]})_{n-1}(z)$ is positive
 definite
with least positive eigenvalue $\epsilon_m$ for all $z\in
 \Omega\setminus
\Omega_{\epsilon_m}$. Thus if $c_m$ is large enough so that
$|\nabla \delta|^2 c_m \ge 2|b_m(z)|$
then, by choosing $\epsilon_m$ small enough, we will  easily see
that the function
$$
\rho_m(z) \equiv  r^{[m+1]}(z)+ c_m \delta(z)^{m+2} \leqno(2.4)
$$
is positive definite in $\Omega\setminus \Omega_{\epsilon_m}$
 and (2.3) holds on 
$\Omega\setminus \Omega_{\epsilon_m}$. Here
$\Omega_t=\{z\in \Omega: \delta(z)<t\}$.
 Then we use arguments in [CKNS]
and [Li1] to extend $\rho_m$
to be defined on $\Omega$ and strictly
pseudoconvex on $\Omega$.
The proof of the theorem is complete. \epf

\medskip

Now we are ready to prove Theorem 1.1.

\medskip

\noindent{\bf Proof of Theorem 1.1.}

Let $ \rho_{4n}$ be the plurisubharmonic defining function
for $\Omega_2$, with $m = 4n$, that we constructed in Theorem 2.2. We let
$$
r(z)= \rho_{4n}(\varphi(z)),\qquad z\in \Omega_1.
$$
Then we have
$$
\det H(r)(z)=\det(H(\rho_{4n})(\varphi(z)))|\det\varphi'(z)|^2=f(z),
\quad z\in \Omega_1
$$
Since $\d \Omega_2 \in C^{\infty}$ 
%Song-Ying:  this is a bit careless.  How exactly are we
%using the fact that $\partial \Omega_2$ is $C^\infty$?  Clearly
%the vanishing of the determinant of the Hessian is making
%things nice for us, but we ought to explain just how.
%This is a critical point!
%Steve: If we can prove the case $\d \Omega$ is $C^2$ then it
%means a lot to us, but we cannot. So it is not important here
% how smooth $\d \Omega$ is. What we really use, may C^{4n}
%or C^{2n}.
then $f(z)^{1/n}\in C^2(\Omegabar)$
 and $f\ge 0$. 
Theorem 2.1 implies that $r\in C^{1,1}(\Omegabar_1)$. So 
$$
\|r\|_{C^{1,1}( \Omegabar_1)}\le C.
$$

It is obvious
that $\exp(\rho_{4n})$ is strictly plurisubharmonic in $\Omega_2$
and that there is a constant $\epsilon > 0$ so that
$$
H(\exp(\rho_{4n})(w) \ge \epsilon I_n,
 \quad \hbox{\rm for\ all }\ \ w\in \Omega_2 .
$$
Thus
\begin{eqnarray*}
{\d^2 \exp( r(z))\over 
\d z_\ell\d \zbar_\ell}
&=&
\sum_{p,\ell=1}^n {\d^2 \exp(\rho_{4n})\over 
\d w_p\d \wbar_{\ell}}(\varphi(z))
{\d \varphi_p\over \d z_\ell}
\overline{ {\d \varphi_{\ell}\over \d z_\ell}}(z)\\
&\ge& \epsilon \sum_{p=1}^n \Big| {\d \varphi_p\over \d z_\ell}(z)\Big|^2.
\end{eqnarray*}
This shows that
$$
\|\varphi\|^2_{{\rm Lip}(\Omegabar_1)} \le 
C\lambda_0^{-1}\|\exp(v)\|_{C^{1,1}(\Omegabar_1)}.
$$
 Thus we have $\det (\varphi'(z))\in L^{\infty}(\Omega_1)$.

If we apply ${\d^2\over \d z_\ell \d z_m}$ to $ r(z) $ 
and use the above result, then we have
$$
\left |\sum_{m=1}^n {\d \rho_{4n}\over \d w_p}(\varphi(z))
{\d^2 \varphi_p\over \d z_\ell \d z_m} \right |
\le C.
$$

Let $z^0\in \Omega_1$ be sufficiently near to $\d \Omega_1$. 
Without loss of
generality, by applying a rotation, 
we may assume that  $z_1^0,\cdots, z_{n-1}^0$ are complex tangential
at $z^0$ and also that at the point $\varphi(z_0)$ the directions 
$ w_1,\cdots, w_{n-1}$ are complex tangential. Thus
\begin{eqnarray*}
\lefteqn{{\d \over \d z_p} \log( \det(\varphi'(z_0)))(z)}\\
&=& \sum_{\ell, m=1}^n \varphi^{\ell m}
 {\d ^2\varphi_m \over \d z_\ell \d z_p}(z)\\
&=& \sum_{\ell<n}\varphi^{\ell m} {\d ^2\varphi_m \over \d z_\ell \d
 z_p}(z)
   + \sum_{m < n}\varphi^{n m} {\d ^2\varphi_m \over \d z_n \d z_p}(z)+
\varphi^{n n} {\d ^2\varphi_n \over \d z_n \d z_p}(z) ,\\ 
\end{eqnarray*}
where $(\varphi^{\ell m})$ is the inverse matrix of $ \varphi'(z)$. 
Since $ \varphi \in \hbox{\rm Lip}_1(\Omegabar_1)$, we have
$$
\left | {\d^2  \varphi_m\over \d z_\ell \d z_p}(z_0)\right |
 \le C \delta_1(z) ^{-1/2}
$$
for all $1\le \ell\le n-1$ and $1\le p\le n$. By (3.8), we have
$$
\left | {\d^2 \varphi_n\over \d z_n \d z_p}(z_0)\right | \le C.
$$

Now we consider the terms with $1\le m\le n-1$. Since
 $(0,\cdots, 0, \varphi_n )$ is normal at 
$\varphi(z_0)$ and $z_1,\cdots, z_{n-1}$ are
complex tangential at $z_0$ to $\d \Omega_1$, we have
$$
\Big| {\d \varphi_n\over \d z_\ell}(z_0)\Big| \le
 C\delta_2(\varphi(z_0))^{1/2}
\le C\delta_1 (z_0)^{1/2}. 
$$ 
Since $\det(\varphi'(z))$ is bounded and 
$$
|\det(\varphi(z)) \varphi^{\ell m}(z)|\le C,
$$
we see that
\begin{eqnarray*} 
\lefteqn{|\det(\varphi'(z)\varphi^{n m}(z_0) {\d ^2 \varphi_m 
\over \d z_n \d z_n}(z_0)|}\\
&\le& C \delta_1(z_0)^{-1} \left |\det ({\d \varphi_p\over \d z_\ell})
_{1\le \ell\le n-1, p\ne m})\right | \\
&\le& C \delta_1(z)^{-1/2}.
\end{eqnarray*}

Combining all the estimates, we have proved that 
$$
|\nabla \det(\varphi'(z))|
=|\det(\varphi'(z) )\nabla\det(\varphi'(z)\log \det (\varphi'(z))|
\le C\delta(z)^{-1/2},
$$
for all $ z\in \Omega_1.$
Hence $\det(\varphi'(z))\in \hbox{\rm Lip}_{1/2}(\Omegabar_1)$.
The proof of Theorem 1.1 is complete. \epf

\section{Reduction of the Mapping Problem}

In this section, we shall state several theorems (from  F.\ Forstneric
[For] and Webster [Web]) which show how to connect the mapping
problem with the reflection principle. 

Let $\Sigma\subset \CC^n$ be a maximal totally real submanifold.  The 
``edge-of-the-wedge'' domain $\Sigma$ is locally defined as follows: 
Near $p\in \Sigma$ we
find $n$ smooth real-valued
functions $r_1,\cdots, r_n$ so that
$$
\Sigma=\{ z: r_1(z)=\cdots =r_n(z)=0\}
$$
and so that the complex gradients
$$
\d r_m=\sum_{p=1}^n {\d r_m\over \d z_p} d z_p
$$
are $\CC$-linearly independent on $\Sigma$.
If $U$ is a neighborhood of $p$ in $\CC^n$, and if 
$\Gamma\subset \RR^n$ is an open convex cone
with vertex zero, then we define the {\it wedge with edge} $\Sigma$ to be
$$
W=W(U, \Gamma)=\{z\in U: r(z)\in \Gamma \}.
$$
Here we have used the notation $r(z) = \bigl ( r_1(z), \dots, r_n(z) \bigr)$.
Then the following theorem is due to Pinchuk/Hasanov [PH] and Coupet [Cou].
Note that, in the statement of the 
theorem, the phrase ``$g$ is asymptotically
holomorphic at $\Sigma$'' has the standard meaning 
that $\overline{\partial} g$
vanishes to infinite order at $\Sigma$.

\begin{theorem}\sl  Let $W=W(U, \Gamma)$ be a wedge with a smooth maximal
totally real edge $\Sigma\subset \CC^n$, let $\Sigma'\subset \CC^{n'}$ be 
another smooth totally real submanifold, and let
$F: W\cup \Sigma\to \CC^{n'}$ be a continuous mapping 
that is smooth on $W$, 
asymptotically holomorphic at $\Sigma$, and such that
$F(\Sigma)\subset \Sigma'$. Then the
restriction of $F$ to $\Sigma \cap U$ is smooth.
\end{theorem}

In order to connect this last theorem to the mapping problem, we need
to introduce an important theorem due to Webster [Web].  First we need some
notation.

Let $M$ and $M'$ be local strictly convex smooth hypersurfaces containing
the origin in $\CC^n$ ($n > 1$). Let $D$ be a domain in $\CC^n$,
with boundary passing through the origin, that is smoothly bounded
and pseudoconvex along $M \subseteq \partial \Omega$, 
and let $D'$ be a similar such domain bounded in part by $M'$.
We assume that we have a biholomorphic mapping $f: D\to D'$ that 
extends to a diffeomorphism on
 $D\cup M$ near the origin and such that this extended
map sends $M$ to $M'$. Since the problem is a local one, we may shrink
our sets toward the origin, and we shall freely do so in the sequel.

Let $a\in \CC^n\setminus \{0\}$ and let
$$
[a]=\{z\in \CC^n :\sum_{m=1}^n z_m a_m=0\}\in \CC\PP^{n-1} .
$$
If $f$ is a holomorphic mapping, then we let
$$
F(z, [a])=(f(z), [a f'(z)^{-1}])
$$
where $f'(z)$ is the Jacobian matrix of $f$ at $z$ and $a f'(z)^{-1}$
is the matrix product of a row vector $a$ with the
 matrix $\bigl [ f'(z)^{-1} \bigr ] $.

Recall that $H_z M \subset T_z M$ is the maximal complex subspace of 
the~$(2n - 1)$-dimensional real 
tangent space $T_z M$. We associate 
to $M$ the smooth submanifold $\tilde{M}$ of 
$X=\CC^n \times \CC P^{n-1}$ defined by
$$
\tilde{M}=\{ (z, H_z M)\in X: z\in M\},
$$
and
$$
\tilde{M'}=\{ (z, H_z M')\in X: z\in M'\}.
$$
The following theorem is due to Webster [Web].

\begin{theorem} \sl If $M\subset \CC^n$ is a strictly pseudoconvex
hypersurface,
then the associated
manifold $\tilde{M}\subset \CC^n\times \CC P^{n-1}$ is totally real
(whence maximally real).
\end{theorem}

If $f$ and $\det (f'(z))$ are  continuous  
up to $D\cup M$, then for each $z\in M$
the derivative $f'(z)$ maps $H_z M$ isomorphically onto $H_{f(z)} M'$.
 Hence the associated mapping $F$ extends continuously from the
domain $\tilde{D}$ to  $\tilde{D} \cup \tilde{M}$ 
and maps $\tilde{M}$ to $\tilde{M'}$.
Its restriction to $\tilde{M}$ is given by
$$
\tilde{f}(z, H_z M)=(f(z), H_{f(z)} M')
$$

When $M=\d \Omega_1$ and $\d \Omega_2$ are strictly pseudoconvex
and smooth then, by Theorem 1.1, we have that $f(z)$ and $\det (f'(z))$
are continuous on $\Omegabar_1$. 
Combining this result with Theorem 3.2 and then with
and 3.1, we obtain a new proof of
the following theorem of C. Fefferman:

\begin{theorem}[C.\ Fefferman]  \sl
Let $\Omega_1$, $\Omega_2$ be strongly pseudoconvex domains
in $\CC^n$ with $C^\infty$ boundaries.  Assume that
$\Phi: \Omega_1 \rightarrow \Omega_2$ is a biholomorphic
(proper) mapping.  Then $\Phi$ extends to a univalent, $C^\infty$
diffeomorphism of $\overline{\Omega}_1$ to $\overline{\Omega}_2$.
\end{theorem}

\noindent {\bf Mailing Address:} 
\vspace*{.1in}

\noindent Department of Mathematics, Washington University,
 St.\ Louis,  MO 63130.
{\sl E-mail:}  sk@math.wustl.edu

\medskip

\noindent
Department of Mathematics, University of California,
Irvine,  CA 92697. \hfill \break
{\sl E-mail:}  sli@math.uci.edu

\end{document}